\documentclass{amsart}
\newtheorem{theorem}{Theorem}
\newtheorem*{theorem*}{Theorem}

\newtheorem{proposition}[theorem]{Proposition}
\newtheorem{corollary}[theorem]{Corollary}
\newtheorem{conjecture}[theorem]{Conjecture}

\theoremstyle{remark}
\newtheorem*{remark*}{Remark}

\theoremstyle{example}
\newtheorem*{example*}{Example}

\begin{document}
\title[Unique Range Sets over Number Fields]
{Yi's Unique Range Set Construction in the Number Field Case}
\author{William Cherry}
\thanks{This article was written while the author was visiting MSRI
while on faculty development leave from UNT.  He would like
to thank both institutions for generous financial support.}
\address{Department of Mathematics\\University of North Texas\\ 
1155 Union Circle \#311430, Denton, TX  76203\\USA}
\email{wcherry@unt.edu}
\date{April 17, 2006 (with minor revisons February 16, 2013)}
\begin{abstract}
H.~X.~Yi's construction of unique range sets for entire functions
is translated to the number theory setting to illustrate that his
construction would work in the number theory setting if one knew
a version of Schmidt's Subspace Theorem with truncated counting functions.
\end{abstract}
\maketitle
\section{Introduction}
A finite set $\mathcal{S}$ of the complex numbers is called a
\textit{unique range set} (counting multiplicity)
for entire functions if whenever $f^*\mathcal{S}=g^*\mathcal{S}$
(pull-back of $\mathcal{S}$ as a divisor)
for two non-constant entire functions $f$ and $g,$ then one must have 
$f=g.$ It is easy to see that a unique range set for entire functions
must contain at least five points.  Indeed, let 
\hbox{$\mathcal{S}=\{a_1,b_1,a_2,b_2\}$}
be a four point set, let $L$ be the M\"obius involution
such that $L(a_j)=b_j$ and $L(b_j)=a_j.$ Then, if $f$ is an entire
function omiting $L(\infty),$ then $L\circ f$ is also an entire function and
\hbox{$(L\circ f)^*\mathcal{S}=f^*(L^*S)=f^*\mathcal{S}.$}
In \cite{Yi}, 
Yi constructed examples of unique range sets for
entire functions of cardinality $\ge 7.$  It seems to be a difficult
open problem to determine if there can be a unique range set 
for entire functions with five or six elements.
\par\smallskip
The main tool in Yi's construction is Nevanlinna theory.  As is now standard,
one can try to transpose Yi's result to number theory.  All constructions
of unique range sets known to me make use
of the Second Main Theorem with truncated counting functions.  Thus,
no existing construction of unique range sets for entire functions will
give an analogous theorem in number theory, except in the case of
number fields with finite unit group.  Some constructions also make
more complicated use of differentiation, which poses an addition challenge
in adapting them to the number theory setting.  However, Yi's 
construction only makes use of the truncated second main theorem, and
thus one can translate his construction to number theory, assuming
a conjecture that Schmidt's Subspace Theorem remains true with
appropriately truncated counting functions.
\par\smallskip
The purpose of this note is to translate Yi's construction to the number
theory setting and to highlight where truncated counting functions are used.
This note contains no new ideas; it is simply a translation of Yi's
paper into number theory.
\par\smallskip
I begin by defining the number theory analogs.  Let $k$ be a number 
field, let $S$ be a finite set of places of $k,$ let $\mathcal{O}_S$
denote the ring of $S$-integers in $k,$ and let $h$ denote an
additive height function on $k.$  Let $\mathcal{S}$ be
a finite set of $S$-units in $k.$  I will call a sequence $x_j$
of $S$-integers in $k$ \textit{admissible} if $h(x_j)\to\infty.$
Two admissible sequences $x_j$ and $y_j$ of $S$-integers
are said to \textit{share}
$\mathcal{S}$ (counting multiplicity)
if for all places $v$ of $k$ not in $S,$ we have
$$
	\left |\prod_{s\in\mathcal{S}}x_j-s\right|_v =
	\left |\prod_{s\in\mathcal{S}}y_j-s\right|_v.
$$
The set $\mathcal{S}$ is called a 
\textit{unique range set} (counting multiplicity) for 
$\mathcal{O}_S$
if whenever $x_j$ and $y_j$ are admissible sequences of $S$-integers
that share $\mathcal{S},$ then $x_j=y_j$ for all but finitely many $j.$
\par\smallskip
From the definition of sharing $\mathcal{S},$ one sees immediately
that it is useful to consider the polynomial
$$
	P_{\mathcal{S}}(X)=\prod_{s\in\mathcal{S}}(X-s).
$$
If $x_j$ and $y_j$ are admissible sequences sharing $\mathcal{S},$
then this precisely means that 
\begin{equation}\label{Peqn}
	P_{\mathcal{S}}(x_j)=u_jP_{\mathcal{S}}(y_j)
\end{equation}
for a sequence
of $S$-units $u_j.$  This leads to the notion of strong uniqueness polynomials.
A polynomial $P$ is called a \textit{strong uniqueness polynomial}
for $\mathcal{O}_S$
if whenever one has two admissible sequences $x_j$ and $y_j$
of $S$-integers such that \hbox{$P(x_j)=cP(y_j)$} for some 
$S$-unit $c,$ then one must have $x_j=y_j$ for all but finitely many $j.$
By Faltings's theorem, one sees that a polynomial $P$ is a strong uniqueness
polynomial for the rings of $S$-integers in all number fields if an only if
the $2$-variable polynomials $P(X)-cP(Y)$ for $c\ne0$
do not have any linear or quadratic factors (over $\mathbf{Q}^a[X,Y]$),
except for the linear factor $X-Y$ when $c=1.$  See \cite{Bilu}, 
\cite{CherryWang}, and \cite{Fujimoto} for various criteria that can
therefore be used to give concrete examples of uniqueness polynomials
for $\mathcal{O}_S.$  Let me also remark here that if the group of
units in $\mathcal{O}_S$ is a finite group, then there is no difference
between the concept of unique range set and strong uniqueness polynomial
because there are only finitely many possibilities of $u_j$
in equation~$(\ref{Peqn}).$
\par\smallskip
As in \cite{Vojta}, for $x$ an element of $k,$
we define the \textit{counting function} (of zeros) by
$$
	N(x)=\frac{1}{[k\colon\mathbf{Q}]}\sum_{v\not\in S}
	\max\{0,\mathrm{ord}_v(x)\}
	[k_v\colon\mathbf{Q}_v]\log p_v,
$$
where $k_v$  denotes the completion of $k$ at the place $v$
and $p_v$ is the prime in $\mathbf{Q}$ which $v$ lies above.
Similarly, the \textit{counting function truncated to multiplicity $\ell,$}
where $\ell$ is a positive integer,
is defined by
$$
	N^{(\ell)}(x)=\frac{1}{[k\colon\mathbf{Q}]}\sum_{v\not\in S}
	\min\{\max\{0,\mathrm{ord}_v(x)\},\ell\}
	[k_v\colon\mathbf{Q}_v]\log p_v,
$$

For Yi's construction to work, one must assume the following 
conjectural strengthening of Schmidt's subspace theorem.
\begin{conjecture}\label{schmidtconj}
Let $L_1,\dots,L_q$ be linear forms in $r+1$-variables with
coefficients in $k$ determining $q$ hyperplanes in
general position in $\mathbf{P}^n.$
Let $x^0_j,\dots,x^r_j$ be $r+1$ sequences of $S$-integers,
at least one of which is admissible, 
such that for each $j$ and each place $v$ not in $S,$
$$
	\max\{|x^0_j|_v,\dots,|x^r_j|_v\}=1,
$$
and such that there is no linear form $L$ such
that $L(x^0_j,\dots,x^r_j)=0$ for infinitely many $j.$
Let $\varepsilon>0.$  Then, for all $j$ sufficiently large,
$$
	(q-r-1-\varepsilon)\max\{h(x^0_j),\dots,h(x^n_j)\}
	\le \sum_{i=1}^q N^{(r)}(L_i(x^0_j,\dots,x^r_j)).
$$
\end{conjecture}

\begin{remark*} If the counting functions on the right were not truncated,
this would be Schmidt's Subspace Theorem. We will only need the conjecture
when $r\le2.$
\end{remark*}

\begin{corollary}\label{schmidtcor}
Assuming Conjecture~\ref{schmidtconj} when $r=1,$ if $x_j$ 
and $y_j$
are sequences of $S$-integers with $x_j$ admissible, and if $A,$ $B,$
and $C$ are non-zero constants such that
$$
	Ax_j+By_j=C
$$
for all but finitely many $j,$ then
$$
	(1-\varepsilon)h(x_j)\le N^{(1)}(x_j)+N^{(1)}(y_j)
$$
for all but finitely many $j.$
\end{corollary}

\begin{proof}
Apply the conjecture with $r=1$ with $x^0_j=1,$ with $x^1_j=x_j,$
and with the three linear forms:
$$
	L_1(x^0,x^1)=x^0,\qquad
	L_2(x^0,x^1)=x^1,\qquad\textnormal{and}\qquad
	L_3(x^0,x^1)=Cx^0-Ax^1.\qedhere
$$
\end{proof}

\section{Yi's Construction}
\begin{theorem}[H.~X.~Yi~{\cite[Theorem~1]{Yi}}]
\label{YiThm}Assume Conjecture~\ref{schmidtconj}
holds when $r\le2.$
Let $n$ and $m$ be relatively prime
positive integers such that $n>2m+4.$ Let $a$ and $b$ be $S$-units
such that the polynomial \hbox{$P(X)=X^n+aX^{n-m}+b$} has no multiple roots
and that the roots of $P$ are $S$-units. Then, the set of zeros 
of $P$ is a unique range set for $\mathcal{O}_S.$
\end{theorem}

\begin{example*} Let $P(X)=X^7+X^6+1.$  Then, assuming
Conjecture~\ref{schmidtconj}, the zeros of $P$
form a unique range set for $\mathcal{O}_S$ for any number field $k$
containing the roots of $P$ and for any finite set of places (containing
all the Archimedean places) $S$ large enough that all the roots of $P$
are $S$-integers.
\end{example*}

\begin{proof}[Proof of Theorem~\ref{YiThm}]
We adopt the convention that throughout the proof all height inequalities
hold for all but finitely many terms and $\varepsilon$ is a positive
number that is adjusted as necessary.
\par\smallskip
Let $\mathcal{S}=\{s_1,\dots,s_n\}$ be the zeros of $P.$
Assume that $x_j$ and $y_j$ are two admissible sequences of
$S$-integers that share $\mathcal{S},$ and so there are $S$-units $u_j$
such that \hbox{$P(x_j)=u_jP(y_j).$}
\par\smallskip
By Roth's Theorem, the
Product formula, and the assumption that $x_j$ and $y_j$ share
$\mathcal{S},$
$$
	(n-1-\varepsilon)h(y_j)\le\sum_{s\in\mathcal{S}}
	N(y_j-s)=\sum_{s\in\mathcal{S}}N(x_j-s)
	\le n h(x_j)+O(1).
$$
Thus by symmetry, we have $h(x_j)=O(h(y_j))$ and $h(y_j)=O(h(x_j)).$
The comparability in height is an important feature of sequences
sharing finite sets.   
Because $u_j=P(x_j)/P(y_j),$ we have by elementary properties of heights,
$$
	h(u_j) \;\le\; h(P(x_j))+h(P(y_j)) = O(h(x_j)).
$$
\par\smallskip
Consider the following auxiliary sequences
\begin{eqnarray*}
	\eta_j &=& -\frac{1}{b}x_j^{n-m}(x_j^m+a)\\
	\zeta_j &=& \frac{1}{b}y_j^{n-m}(y_j^m+a)u_j
\end{eqnarray*}
Because $n$ is somewhat larger than $m$ and $\eta_j$ and $\zeta_j$
begin with something to the $n-m$ power, they will have places
dividing them with moderate multiplicity.  Exploiting this 
extra multiplicity
is where we will use Conjecture~\ref{schmidtconj}, and finding a way
to exploit this sort of multiplicity without referring to an unproven
conjecture is the main obstacle in using existing analytic 
constructions of unique range sets in the number field setting.
\par\smallskip
Now, notice that
$$
	\eta_j+u_j+\zeta_j = \frac{-1}{b}(P(x_j)-b)+u_j
	+\frac{1}{b}(P(x_j)-bu_j)=1.
$$
\par\smallskip
Assume for the moment that there is no linear form in three 
variables $L$ such that $L(\eta_j,u_j,\zeta_j)=0$ for 
infinitely many $j.$  Then, we may apply Conjecture~\ref{schmidtconj}
to conclude that
$$
	(1-\varepsilon)\max\{h(\eta_j),h(u_j),h(\zeta_j)\}
	\le N^{(2)}(\eta_j)+
	N^{(2)}(u_j)+
	N^{(2)}(\zeta_j)+N^{(2)}(\eta_j+u_j+\zeta_j).
$$
Because $u_j$ is an $S$-unit and $\eta_j+u_j+\zeta_j=1,$ 
$N^{(2)}(u_j)=N^{(2)}(\eta_j+u_j+\zeta_j)=0.$  Clearly,
$$
	N^{(2)}(\eta_j)\le 2N^{(1)}(x_j)+N(x_j^m+a) \qquad\textnormal{and}
	\qquad
	N^{(2)}(\zeta_j)\le 2N^{(1)}(y_j)+N(y_j^m+a).
$$
Because the counting functions are bounded by the heights (Product Formula)
and $h(x_j^m+a)$ is comparable to $mh(x_j)$ and similarly for $y_j,$
we get
$$
	(1-\varepsilon)\max\{h(\eta_j),h(u_j),h(\zeta_j)\}
	\le (2+m)h(x_j)+(2+m)h(y_j) + O(1).
$$
Also,
$$
	h(\eta_j)\ge n h(x_j)+O(1),
$$
so
$$
	(n-\varepsilon)h(x_j)\le(2+m)h(x_j)+(2+m)h(y_j)+O(1).
$$
Reversing the roles of $x_j$ and $y_j$ we also get
$$
	(n-\varepsilon)h(y_j)\le(2+m)h(x_j)+(2+m)h(y_j)+O(1).
$$
Adding the previous two inequalities give
$$
	(n-\varepsilon)(h(x_j)+h(y_j))\le(4+2m)(h(x_j)+h(y_j))+O(1),
$$
which contradicts the assumption that $n>2m+4.$
\par\smallskip
The proof is completed by the following proposition. 
\renewcommand{\qed}{}
\end{proof}

\begin{proposition}
Let $c_1,$ $c_2$ and $c_3$ be in $k$ and not all zero.  
With the notation as in the theorem, if $x_j\ne y_j$
for infinitely many $j,$ then 
$$
	c_1\eta_j+c_2u_j+c_3\zeta_j=0
$$
for at most finitely many $j.$
\end{proposition}

\begin{proof}
Suppose $c_1\eta_j+c_2u_j+c_3\zeta_j=0$ for infinitely many $j.$
Clearly at least two of the $c_i$ are non-zero.
\par\smallskip
\textit{Case $c_1=0$:} Then, $\zeta_j=(c_2/c_3)u_j$
and so $y_j^{n-m}(y_j^m+a)=bc_2/c_3,$ which contradicts the
assumption that $h(y_j)\to\infty.$ So this case does not occur.
\par\smallskip
We may now assume $c_1\ne0.$  Because $\eta_j+u_j+\zeta_j=1,$ we can
remove $\eta_j$ to get $C_2u_j+C_3\zeta_j=1,$ where
$$
	C_i=1-\frac{c_i}{c_1}.
$$
\par\smallskip
\textit{Case $C_2\ne0$ and $C_3\ne0$:} In this case,
$$
	C_3\zeta_ju_j^{-1}-u_j^{-1}=-C_2,
$$
and so we can apply Corollary~\ref{schmidtcor} to conclude
$$
	(1-\varepsilon)h(\zeta_ju_j^{-1})\le N^{(1)}(\zeta_ju_j^{-1})
$$
for infinitely many $j.$  But,
$$
	\zeta_ju_j^{-1}=\frac{1}{b}y_j^{n-m}(y_j^m+a),
$$
and hence $h(\zeta_ju_j^{-1})=nh(y_j)+O(1).$ Also,
$$
	N^{(1)}(\zeta_ju_j^{-1})\le
	N^{(1)}(y_j)+N^{(1)}(y_j^m+a) \le (m+1)h(y_j)+O(1).
$$
This contradicts $n>2m+4.$
\par\smallskip
\textit{Case $C_2=0$:} In this case,
$$
	y_j^{n-m}(y_j^m+a)=bC_3^{-1}u_j^{-1}.
$$
Enlarging $S$ if necessary, we may assume $bC_3^{-1}$ is a unit in
$\mathcal{O}_S.$  This implies that $y_j$ and $y_j^m+a$ are $S$-units for
all $j,$ which is a contradiction to the $S$-integer version of Picard's
theorem.
\par\smallskip
\textit{Case $C_3=0$:}  In this case we see that $u_j=C_2^{-1}$ for
infinitely many $j.$  Thus, either $x_j=y_j$ for all but finitely
many of these $j,$ or $P$ is not a strong uniqueness polynomial.
We have already remarked that if $P$ is a strong uniqueness polynomial
for entire functions, then it is also a strong uniqueness polynomial
for $S$-integers.  Thus, the proof is completed by showing 
$P$ is a strong uniqueness polynomial as in Yi \cite{Yi}, or
alternatively as in \cite{Fujimoto}.  
This is where the assumption that $n$ and $m$
are relatively prime is used.
\end{proof}

\end{document}